\documentclass{amsart}

\usepackage[latin1]{inputenc}
\usepackage{latexsym}
\usepackage{amssymb}
\usepackage[T1]{fontenc}	
\usepackage[english]{babel}
\usepackage{amsmath,amssymb}
\usepackage{amsfonts}
\usepackage{amscd}
\usepackage{empheq}
\usepackage{fancyvrb}
\usepackage{marvosym}
\usepackage{mathrsfs}
\usepackage{amsthm}
\usepackage{verbatim}
\usepackage{tikz}
\usetikzlibrary{chains}
\usepackage{color}
\usepackage{chngpage}
\usetikzlibrary{matrix,through,arrows}

\renewcommand{\d}{\delta }
\newcommand{\D }{\Delta }

\newcommand{\e }{\varepsilon }

\renewcommand{\l }{\lambda }

\newcommand{\n }{\nabla }

\renewcommand{\o }{\omega }

\newcommand{\be}{\begin{equation}}
\newcommand{\ee}{\end{equation}}
\newenvironment{pf}{\noindent{\sc Proof}.\enspace}{\rule{2mm}{2mm}\medskip}
\newenvironment{pfn}{\noindent{\sc Proof}}{\rule{2mm}{2mm}\medskip}
\newtheorem{thm}{Theorem}[section]
\newtheorem{pro}[thm]{Proposition}
\newtheorem{definition}[thm]{Definition}
\newtheorem{lem}[thm]{Lemma}
\newtheorem{fct}[]{Fact}
\newtheorem{rem}[thm]{Remark}
\newtheorem{cor}[thm]{Corollary}

\newcommand{\Z}{\mathbb{Z}}

\renewcommand{\epsilon}{\varepsilon}

\begin{document}

\title{ON THE SOLVABILITY OF SINGULAR LIOUVILLE EQUATIONS ON COMPACT SURFACES OF ARBITRARY GENUS}
\author{Alessandro Carlotto}
\address{Department of Mathematics, Stanford University, 
94305 Stanford CA}
\email{carlotto@stanford.edu}
\urladdr{www.math.stanford.edu/$\sim$carlotto}



\begin{abstract}
In the first part of this article, we complete the program announced in the preliminary note \cite{cm} by proving a conjecture presented in \cite{cm2} that states the equivalence of contractibility and $p_{1}-$stability for generalized spaces of formal barycenters $\Sigma_{\rho,\underline{\alpha}}$ and hence we get \textsl{purely algebraic} conditions for the solvability of the singular Liouville equation on Riemann surfaces. This relies on a structure decomposition theorem for $\Sigma_{\rho,\underline{\alpha}}$ in terms of maximal strata, and on elementary combinatorial arguments based on the selection rules that define such spaces. Moreover, we also show that these solvability conditions on the parameters are not only \textsl{sufficient}, but also \textsl{necessary} at least when for some $i\in\left\{1,\ldots,m\right\}$ the value $\alpha_{i}$ approaches $-1$. This disproves a conjecture made in Section 3 of \cite{tars} and gives the first non-existence result for this class of PDEs without any genus restriction. The argument we present is based on a combined use of maximum/comparison principle and of a Pohozaev type identity and applies for arbitrary choice both of the underlying metric $g$ of $\Sigma$ and of the datum $h$.
\end{abstract}

\maketitle

\

\section{Introduction and statement of the results}

\medskip

\noindent

Let $\Sigma$ be a closed\footnote{Namely $\Sigma$ will be always assumed to be compact and with no boundary.}, orientable surface, which is endowed with a Riemannian metric $g$ and let $\rho$ be a real parameter: then it is of wide interest the study of the following partial differential equation
\begin{equation}
\label{reg}
-\Delta_{g}u+\rho=\rho\frac{h(x)e^{2u}}{\int_{\Sigma}h(x)e^{2u}\,dV_{g}} \quad \textrm{on}\ \Sigma,
\end{equation}
that are dated back even to Liouville. Mainly motivated by Differential Geometry, namely by the central problem of prescribing the Gaussian curvature under a \textsl{conformal} change of metric $g\mapsto e^{2u}g$, these equations have been extensively studied both by topological and variational methods (see for instance \cite{clin}, \cite{djlw1}, \cite{zind}).\newline
More recently, the \textsl{singular} counterpart of \eqref{reg} has been the object of several works and specifically
in the recent articles \cite{cm2} and \cite{mr} the authors present a variational approach to the problem
\begin{equation}
\label{sing}
-\Delta_{g}u=\rho\left(\frac{h(x)e^{2u}}{\int_{\Sigma}h(x)e^{2u}\,dV_{g}}-1\right)-2\pi\sum_{i=1}^{m}\alpha_{j}(\delta_{p_{j}}-1)
\end{equation}
where $p_{j}\in\Sigma$ are some fixed points, $h$ is a smooth positive function, while $\rho$ and $\alpha_{1},\ldots, \alpha_{m}$ are real numbers. In fact, we will restrict our attention to the case studied in \cite{cm2}, namely $\alpha_{j}\in\left(-1,0\right)$ for all $j=1,2,\ldots,m$. Moreover, notice that a necessary condition for the solvability of both \eqref{reg} and \eqref{sing} (as well as of its \textsl{regularised} version \eqref{mod}) is that $Vol_{g}\left(\Sigma\right)=1$, which will be always implicitly assumed in the sequel.

This equation also has a strong geometric flavour, since it extends the question presented above to the case when the \textsl{target} curvature is given by the sum of a smooth function and a finite linear combination of Dirac masses describing conical singularities. Indeed, it can be shown that that finding a solution for \eqref{sing} (with an appropriate choice of $\rho$ and $h$) is the same as prescribing the Gaussian curvature on a surface with conical singularities. The details of this formulation of the problem have been discussed for instance in Proposition 1.2 in \cite{bdm}. Here we only recall that a cone of angle $\phi=2\pi(1+\alpha)$ at a point $p$ is described in \eqref{sing} by a summand $-2\pi\alpha\delta_{p}$, as can justified by means of a well-known extension of the Gauss-Bonnet formula (see \cite{tr}). 

Equation \eqref{sing} also arises in the study of self-dual multivortices in the Electroweak Theory by Glashow-Salam-Weinberg, where $u$ can be interpreted as the logarithm of the absolute value of the wave function and the points $p_{j}$'s are the \textsl{vortices}, where the wave function vanishes. This class of problems has proved to be relevant in other physical frameworks, such as the study of the statistical mechanics of point vortices in the mean field limit and the abelian Chern-Simons Theory. An excellent account on the different perspectives of the recent research concerning Liouville type equations on surfaces is given in \cite{tars} (see also the ample references therein).

One basic principle in studying equation \eqref{sing} is trying to exploit its variational structure. Indeed, if we perform the \textsl{singular} change of variable $\widetilde{u}:=u-2\pi\sum_{j=1}^{m}\alpha_{j}G_{p_{j}}$ (where $G_{p_{\ast}}\left(\cdot\right)=G\left(\cdot,p_{\ast}\right)$ are suitable Green functions, see Proposition \ref{verde}), then \eqref{sing} transforms into
\begin{equation}\label{mod}
-\Delta_{g}\widetilde{u}+\rho=\rho\frac{\widetilde{h}(x)e^{2\widetilde{u}}}{\int_{\Sigma}\widetilde{h}(x)e^{2\widetilde{u}}\,dV_{g}} \quad \textrm{on}\ \Sigma,
\end{equation}
with $\widetilde{h}(x)=h(x)e^{4\pi\sum_{j=1}^{m}\alpha_{j}G_{p_{j}}}$ and, by comparison with the \textsl{regular} case (see \cite{zind}), we get that \eqref{mod} is nothing but the Euler-Lagrange equation for the modified functional
\begin{equation}
\label{functm}
J_{\rho,\underline{\alpha}}(\widetilde{u})=\int_{\Sigma}\left|\nabla_{g}\widetilde{u}\right|^{2}\,dV_{g}+2\rho\int_{\Sigma}\widetilde{u}\,dV_{g}
-\rho\log\int_{\Sigma}\widetilde{h}(x)e^{2\widetilde{u}}\,dV_{g}, \qquad \widetilde{u}\in H^{1}(\Sigma,g)
\end{equation}
(where $\underline{\alpha}=\left(\alpha_{1},\ldots,\alpha_{m}\right)\in\left(-1,0\right)^{m}$). In fact, there are no obstructions in reducing the solution of $\eqref{sing}$ to the study of $\eqref{mod}$, the regularity issues being described in Section 5 of \cite{cm2}.

The key idea, to that aim, is trying to detect a change in the topology of the sublevels $J_{\rho,\underline{\alpha}}^{L}$ (see the \textsl{Notations} paragraph at the end of this section) of our functional \eqref{functm}. Specifically, it has been shown that the \textsl{very high} sublevels are a deformation retract of the whole ambient space $H^{1}\left(\Sigma, g\right)$ while the \textsl{very low} sublevels are described (in the sense of homotopy equivalence) by a rather non-trivial class of spaces, that might be called \textsl{generalized spaces of formal barycenters}. More precisely, it was proved in \cite{cm2} (using some techniques introduced in \cite{dm}) that such (partial) homotopy equivalence \footnote{See Lemma 4.3 in \cite{cm2} for a precise statement.} exists whenever the parameters $\alpha_{1},\ldots,\alpha_{m}$ belong to the real interval $\left(-1,0\right)$ and $\rho$ is not a singular value for \eqref{sing} in the following sense.

\begin{definition}
We say that $\overline{\rho}>0$ is a singular value for Problem \eqref{sing} if
\begin{equation}
\label{defsv}
\overline{\rho}=4\pi n+4\pi\sum_{i\in I}(1+\alpha_{i})
\end{equation}
for some $n\in\mathbb{N}$ and $I\subseteq\left\{1,\ldots,m\right\}$ (possibly empty) satisfying
$n+card\left(I\right)>0$. The set of singular values will be denoted by $\mathfrak{S}=\mathfrak{S}\left(\underline{\alpha}\right)$.
\end{definition}

Let us see then \textsl{how} such model spaces $\Sigma_{\rho,\underline{\alpha}}$ are actually defined.

\begin{definition}
Given a point $q\in\Sigma$ we define its weighted cardinality as follows:
\begin{displaymath}
\chi(q)=\left\{\begin{matrix} 1+\alpha_{i} & \mbox{if }q=p_{i} \ \mbox{for some} \ i=1,\ldots,m; \\ 1 & \mbox{otherwise}.
\end{matrix}\right.
\end{displaymath}
The cardinality of any finite set of (pairwise distinct) points on $\Sigma$ is obtained extending $\chi$ by additivity.
\end{definition}
This enables us to easily describe selection rules to determine admissibility conditions for specific barycentric configurations in dependence on the values of the $\alpha_{i}$'s and $\rho.$

\begin{definition}
\label{selru}
Suppose all the parameters $\rho,\alpha_{1},\ldots,\alpha_{m}$ are fixed. We define the corresponding space of formal barycenters as follows
\begin{equation}
\Sigma_{\rho,\underline{\alpha}}=\left\{\sum_{q_{j}\in J}t_{j}\delta_{q_{j}}:\ J \ \textrm{is finite,} \ \sum_{q_{j}\in J}t_{j}=1,\ t_{j}\geq0,\ q_{j}\in\Sigma, \quad 4\pi\chi(J)<\rho\right\}.
\end{equation}
\end{definition}

We will consider $\Sigma_{\rho,\underline{\alpha}}$ endowed with the weak topology corresponding to the duality with $C^{1}(\Sigma,g)$. It is easy to see that such topology is equivalently determined by the distance function
\begin{displaymath}
d:\Sigma_{\rho,\underline{\alpha}}\times \Sigma_{\rho,\underline{\alpha}} \to \mathbb{R}_{\geq 0}\ , \quad d(\sigma_{1},\sigma_{2})= \sup_{\left\|f\right\|_{C^{1}\left(\Sigma\right)}\leq 1}\left(\sigma_{1}-\sigma_{2},f\right).
\end{displaymath}
This will be a useful tool to perform some explicit computations. Notice that since we are considering negative weights the topological structure of $\Sigma_{\rho,\underline{\alpha}}$
is in general richer than that of 
\begin{displaymath}
\Sigma^{k,\emptyset}=\left\{\sum_{j=1}^{k}t_{j}\delta_{q_{j}}:\  \sum_{j=1}^{k}t_{j}=1,\ t_{j}\geq0,\ q_{j}\in\Sigma\right\},
\end{displaymath}
(which are the corresponding model spaces for the regular problem) and strongly depends on the values of the parameters $\rho$ and $\underline{\alpha}$. For instance, when $m=2, \alpha_{1}=\alpha_{2}=\alpha$, $\rho>8\pi\left(1+\alpha\right)$ and $ 4\pi<\rho<4\pi(2+\alpha)$ we have that $\Sigma_{\rho, \underline{\alpha}}$ can be visualised as a space obtained by gluing together a mirror image of $\Sigma$ and a linear handle joining the singular points $p_{1}$ and $p_{2}$.

Due to all the previous remarks it should be clear that \textsl{the key issue here is determining whether, for given values of $\rho$ and $\underline{\alpha}$ the space $\Sigma_{\rho,\underline{\alpha}}$ is contractible or not}. Moving from the study of a wide range of examples, it was conjectured in \cite{cm2} that contractibility of such a space is in fact equivalent to its $p_{1}-$stability, which is defined as follows. 

\begin{definition}
\label{piuno}
Given the parameters $\rho$ and $\underline{\alpha}$, we say that the corresponding model space $\Sigma_{\rho,\underline{\alpha}}$ is $p_{j}-$stable for some index $j\in{1,2,\ldots,m}$ if one of the following two equivalent conditions holds:
\begin{enumerate}
\item{Whenever $\sigma\in \Sigma_{\rho,\underline{\alpha}}$ then $(1-t)\sigma+t\delta_{p_{j}} \in \Sigma_{\rho,\underline{\alpha}} \ \forall\ t \in \left[0,1\right]$;}
\item{Whenever $k\in\mathbb{N}$ and a set $I$ are such that \begin{displaymath}4\pi\left[k+\sum_{i\in I}\left(1+\alpha_{i}\right)\right]<\rho \end{displaymath} then also \begin{displaymath}4\pi\left[k+\sum_{i\in\left\{j\right\}\cup I}\left(1+\alpha_{i}\right)\right]<\rho.\end{displaymath}}
\end{enumerate}
\end{definition}

Obviously, if $\Sigma_{\rho,\underline{\alpha}}$ is $p_{1}-$stable then it is contractible or, more precisely, it deformation-retracts onto $\delta_{p_{1}}$ in the ambient space $C^{1}\left(\Sigma,g\right)^{\ast}$ by means of the homotopy map $H:\Sigma_{\rho,\underline{\alpha}}\times\left[0,1\right]\rightarrow \Sigma_{\rho,\underline{\alpha}}$ given by $H(\sigma,t)=(1-t)\sigma+t\delta_{p_{1}}.$ The converse is much less clear and it is the object of the first part of this article.

\begin{thm}[topological version]\label{cong}
The space of formal barycenters $\Sigma_{\rho,\underline{\alpha}}$ is contractible if and only if it is $p_{1}-$stable.
\end{thm}

Such theorem can be immediately turned into algebraic, namely algorithmic, form.

\begin{thm}[algebraic version]
\label{conalg}
The space of formal barycenters $\Sigma_{\rho,\underline{\alpha}}$ is  NOT contractible if and only if there exist a number $k\in\mathbb{N}$ and a set $I\subseteq\left\{1,2,\ldots,m\right\}$, possibly empty, such that $k+card(I)>0$ and 
\begin{displaymath}
\rho>4\pi\left[k+\sum_{i\in I}\left(1+\alpha_{i}\right)\right] \ \wedge \ \rho<4\pi\left[k+\sum_{i\in\left\{1\right\}\cup I}\left(1+\alpha_{i}\right)\right].
\end{displaymath}
\end{thm}

In Section 2 of this article, we list a number of basic properties concerning the spaces $\Sigma_{\rho,\underline{\alpha}}$, mainly dealing with its structure of a stratified space. Using these preliminary results, we prove this theorem in Section 3: this completes the program described and announced in the preliminary note \cite{cm}. Rather surprisingly, as a net result of our theory we get a general existence result for equation \eqref{sing} based on a sufficient condition of purely algebraic nature (with respect to the values of $\rho$ and $\alpha_{1},\ldots, \alpha_{m}$).

\begin{thm}
\label{final}
Whenever the parameters $\underline{\alpha}\in\left(-1,0\right)^{m}$ and $\rho\in\mathbb{R}\setminus\mathfrak{S}\left(\underline{\alpha}\right)$ are chosen so that the space $\Sigma_{\rho,\underline{\alpha}}$ is NOT $p_{1}-$stable, then equation \eqref{sing} is solvable.
\end{thm}

This is an immediate consequence of Theorem \ref{cong} (or, equivalently, \ref{conalg}) and Theorem 1.6 in \cite{cm2}.

\begin{rem}
The regularity of the solutions has been studied in Section 5 of \cite{cm2}: indeed it was proved that when the space $\Sigma_{\rho,\underline{\alpha}}$ is \textsl{not} contractible, then Problem \eqref{sing} admits a solution $u$ such that we can write $u=v+\sum_{j=1}^{m}\alpha_{j}G_{p_{j}}$ where $G_{p_{\ast}}$ refers to the Green functions introduced above and $v\in C^{\gamma}(\Sigma,g)$, for any $\gamma\in\left[0,\gamma_{0}\right)$ with $\gamma_{0}\in\left(0,1\right)$, solving equation \eqref{mod}. In fact, such $u$ will be smooth ($C^{\infty}$) away from the singularities and will have logarithmic blow-ups at the points $p_{j}$ for $j=1,\ldots,m$ since (see \cite{aub}) $G_{p_{j}}\left(x\right)\simeq \log d_{g}\left(x,p_{j}\right)$ for $x\to p_{j}$.
\end{rem}

The peculiar character of this result is reducing an \textsl{analytical} problem to an \textsl{algebraic} criterion, which in turns is derived by means of a topological approach: what is proved in this article allows to completely by-pass explicit topological assumptions in this final statement. \newline

A rather remarkable consequence of our main result concerning the equivalence between contractibility and $p_{1}-$stability of $\Sigma_{\rho,\underline{\alpha}}$ is the following fact, which seems highly non-trivial from the purely topological viewpoint.

\begin{cor}
\label{ascontr}
For any index $j\in{1,\ldots,m}$ and value of $\rho\in\mathbb{R}\setminus \mathfrak{S}\left(\alpha_{1},\ldots,\alpha_{j-1},\alpha_{j+1},\ldots,\alpha_{m}\right)$,
there exists a threshold value $\alpha^{\ast}$ (depending on $\rho$) such that the corresponding space $\Sigma_{\rho,\underline{\alpha}}$ is contractible for $\alpha_{j}<\alpha^{\ast}$. 
Here we have assumed, by extension of the previous definition given in \eqref{defsv} that
\begin{equation}\label{sve}
 \mathfrak{S}\left(\alpha_{1},\ldots,\alpha_{j-1},\alpha_{j+1},\ldots,\alpha_{m}\right)=\left\{4\pi\left(n+\sum_{I}(1+\alpha_{i})\right) \textrm{for} \ n\in \mathbb{N} \ \textrm{and} \  I\subseteq\left\{1,\dots,j-1,j+1,\ldots,m\right\}\right\}.
\end{equation}
\end{cor}

\medskip

As a logical counterpart of the previous discussion, we might naturally ask what can in fact be said in case the values of $\rho$ and $\underline{\alpha}$ are chosen so that $\Sigma_{\rho,\underline{\alpha}}$ is contractible. Now, the reader should understand that it is quite unreasonable to ask for the previous algebraic conditions to be also necessary for solvability of \eqref{sing}. Indeed, the min-max approach which lies at the basis of our variational theory only detects saddle points corresponding to changes in the topology of the sublevels of $J_{\rho,\underline{\alpha}}$ while for instance it neglects \textsl{local minima} whose existence can be obtained, at least in some special cases, even by an appropriate smart choice of the datum $h$ (see e.g. \cite{dmfr}).  Nevertheless, it is certainly an interesting question to ask whether at least some partial converse of Theorem \ref{final} holds.
 
The main purpose of the second part of this article is to try to answer such a question.  Moving from a recent conjecture by G. Tarantello (\cite{tars}), we first study the case $m=1$, namely let $p\in\Sigma$ be the only singular point for \eqref{sing} and let $\alpha\in\left(-1,0\right)$ be the associated angular parameter. For $\rho\in\left(0,4\pi\left(1+\alpha\right)\right)$ the functional $J_{\rho,\alpha}$ is coercive (thanks to an inequality due to Troyanov, see \cite{tr}), while for $\rho\in\left(4\pi\left(1+\alpha\right),4\pi\right)$ the space $\Sigma_{\rho.\alpha}$ is indeed contractible and the sovability issues immediately become not trivial at all. In fact, concerning this regime, it was proved in the above mentioned article that when $\Sigma=\mathbb{S}^{2}$ equation \eqref{sing} is NOT solvable. The proof relies on the use of the stereographic projection and of a Pohozaev-type identity. \textsl{In such regime} our theory does \textsl{not} distinguish between $\mathbb{S}^{2}$ and positive genus surfaces and so we were motivated to believe that such theorem could actually be proved without topological restrictions on $\Sigma$. Unfortunately, the method used by Tarantello is highly peculiar and does not extend to the case when the genus of $\Sigma$ is positive and indeed it was conjectured in \cite{tars}, in the context of the degree theory approach to this problem, that \eqref{sing} was indeed solvable on tori for $\rho\in\left(4\pi\left(1+\alpha\right),4\pi\right)$. In Section 4 and 5 of this article, we disprove this fact by showing the following more general result.

\begin{thm}
\label{bh}
Assume that $\alpha_{i}\in\left(-1,0\right)$ for $i\in\left\{1,\ldots,m\right\}$, and let $\rho\notin \mathfrak{S}\left(\alpha_{1},\ldots,\alpha_{j-1},\alpha_{j+1},\ldots,\alpha_{m}\right)$ for the corresponding values of $\alpha_{1},\ldots,\alpha_{j-1},\alpha_{j+1},\ldots,\alpha_{m}$ (see \eqref{sve}). Then there exists $\alpha_{\ast}\in\left(-1,0\right)$ such that equation $\eqref{sing}_{\rho,\underline{\alpha}}$ is not solvable for $\alpha_{j}\leq\alpha_{\ast}$. Moreover, such $\alpha_{\ast}$ can be chosen uniformly for $\rho$ in any relatively compact subset of $\mathbb{R}_{>0}\setminus \mathfrak{S}\left(\alpha_{1},\ldots,\alpha_{j-1},\alpha_{j+1},\ldots,\alpha_{m}\right) $.
\end{thm}

This is proved by contradiction, namely we study the behaviour of a sequence of (normalized) solutions of \eqref{sing} for fixed $\rho$ and $\alpha_{j}$ approaching the threshold value $-1$. Our argument is based on a rather \textsl{soft} use of the classical blow-up analysis by Brezis and Merle (\cite{bm}) away from the singularity (and in fact on its extension given in \cite{bmt}), on a maximum/comparison principle to rule out lower boundedness of any subsequence at (any) positive distance from $p$ and finally on a \textsl{local} version of the Pohozaev identity.  The proof we give has the following distinctive features:
\begin{itemize}
\item{it bypasses an explicit blow-up analysis for equation \eqref{sing} for fixed $\rho$ and $\alpha_{j}\searrow -1$, which is though partially obtained as a byproduct;}
\item{it works without any restriction on the genus of $\Sigma$, so that the case of the sphere is also covered here with no reference to the stereographic projection;}
\item{it requires no assumption on the Riemannian metric $g$ and the datum $h$.}
\end{itemize}

Section 4 is devoted to a detailed proof of this result in the special case when $m=1$ and $\rho\in\left(0,4\pi\right)$, while a description of the extension of our argument to the \textsl{general} situation of Theorem \ref{bh} is given in Section 5.\newline
Finally, let us say that the above results admit a vague but suggestive physical interpretation: indeed when some of the parameters $\alpha_{j}$ approaches $-1$, the corresponding singularity acts as a sort of \textsl{black hole} that causes an energy loss which is formally reflected by the (failure of) the Pohozaev identity.

For the sake of brevity, we do not discuss here the specific consequences of our results with respect to the geometric problem described above (first studied in \cite{kw} for the regular case): this is the object of a forthcoming paper with A. Malchiodi. However, the effectiveness of the variational approach to such problem has been recently shown in \cite{bdm} (under the assumption that the genus of $\Sigma$ is positive and $\alpha_{j}\geq 0$ for $j=1,\ldots,m$), where general existence and multiplicity results are obtained. 

\medskip

\textsl{Notations}. Throughout this article, we will assume that the labeling of the singular points is chosen so that $\alpha_{1}\leq \alpha_{2}\leq\ldots\leq\alpha_{m}.$ The dependence on the parameters $\alpha_{1},\ldots,\alpha_{m}$ of the modified datum $\widetilde{h}$ for problem \eqref{mod} will be implicit, as explained in Remark \ref{indic}. We refer to the sublevels of the functional $J_{\rho,\underline{\alpha}}$ by $J_{\rho,\underline{\alpha}}^{L}$, this being the subset of $H^{1}\left(\Sigma,g\right)$ of the functions $u$ for which $J_{\rho,\underline{\alpha}}\left(u\right)\leq L$. Large positive constants are always denoted by $C$ and the exact value of $C$ is allowed to vary from formula to formula and also within the same line. When we want to stress the dependence on some parameter, we add subscripts to $C$, hence obtaining things like $C_{\eta}$ and so on. Landau convention will also be adopted, so that for instance $o_{\d}(1)$ stands for a function of $k\in\mathbb{N}$ converging to 0 as $k\to+\infty$ and depending on a parameter $\d$.  The Lebesgue measure on $\mathbb{R}^{2}$ will be $\mathscr{L}^{2}$. In applying the Mayer-Vietoris Theorem, we will regularly denote by $A$ and $B$ the sets we work with, so that our space splits as $X=A\cup B$: it is well known that $A$ and $B$ are tipically required to be \textsl{open} subspaces of $X$. Nevertheless, in order to avoid unnecessariliy tedious notation, we will work with $A$ and/or $B$ possibly not open in $X$, in this case tacitly assuming that suitable open neighbourhoods of $A, B$ should be considered instead. Finally, in this note we will always deal with homology taking coefficients in $\Z_{2}:$ as usual in the Algebraic Topology literature, $\mathbb{Z}_{2}$ denotes the field with two elements. In order to avoid ambiguities in some of the definitions and in the statement of our results, we remark that $\mathbb{N}=\left\{0,1,2,\ldots\right\}$.

\medskip

\section{The structure of $\Sigma_{\rho,\underline{\alpha}}$ as a stratified space}

\medskip

\noindent

We need to start by introducing some notation.
For $k \in \mathbb{N}$ and a set of indices $I\subseteq \left\{1,\ldots,m\right\}$ satisfying the relation $4\pi \left[k+\sum_{I}\left(1+\alpha_{i}\right)\right]<\rho$ we define the set
\begin{displaymath}
\Sigma^{k,I}=\left\{\sum_{I}{s_{i}\delta_{p_{i}}}+\sum_{j=1}^{k}t_{j}\delta_{q_{j}}\right\},
\end{displaymath}
where
\begin{itemize}
\item{$s_{i}\in\left[0,1\right]$ for any $i\in I$;}
\item{$t_{j}\in\left[0,1\right]$ for any $j=1,\ldots,k$;}
\item{$\sum_{i}s_{i}+\sum_{j}t_{j}=1$;}
\item{$q_{j}\in \Sigma,$ for any $j=1,\ldots,k$.}
\end{itemize}

\begin{definition}
\label{ordstr}
Given two couples $\left(k_{1},I_{1}\right)$ and $\left(k_{2},I_{2}\right)$, we will write that $\Sigma^{k_{1},I_{1}}\preceq \Sigma^{k_{2},I_{2}}$ in case $\Sigma^{k_{1},I_{1}}\subseteq \Sigma^{k_{2},I_{2}}$ or, equivalently, if $k_{2}\geq k_{1}$ and the set $I_{1}$ can be split into two subsets, say $I_{1}'$ and $I_{1}''$, such that:
\begin{itemize}
\item{$I_{1}'\subseteq I_{2}$;}
\item{$card\left(I_{1}''\right)\leq k_{2}-k_{1}$}.
\end{itemize}
We will denote by $\prec$ the strict order associated to $\preceq$.
\end{definition}

Notice that the \textsl{partial} order $\prec$ restricts to a \textsl{total} order on the class of strata $\Sigma^{k,\emptyset}$, namely
\begin{displaymath}
\Sigma^{1,\emptyset}\prec\Sigma^{2,\emptyset}\prec\ldots\prec\Sigma^{k-1,\emptyset}\prec\Sigma^{k,\emptyset}, \quad \rho\in\left(4k\pi,4\left(k+1\right)\pi\right),
\end{displaymath}
which corresponds to the \textsl{regular} problem, that is in absence of singularities.

This notion of \textsl{partial order} among strata of a given space $\Sigma_{\rho,\underline{\alpha}}$ allows to associate to such space (that is, in practice, to a given choice of $\alpha_{1},\ldots,\alpha_{m}$ and non-singular $\rho$) a combinatorial \textsl{oriented graph} that we will call, for brevity \textsl{S-graph}. Obviously, we mean that each stratum is represented by a node and inclusion is described by an oriented edge. We will say that a stratum in $\Sigma_{\rho,\underline{\alpha}}$ is \textsl{maximal} if it is maximal with respect to the order $\prec$ on the family of strata of $\Sigma_{\rho,\underline{\alpha}}$.  We will make use, in the sequel, of the following easy but crucial result.

\begin{fct}\label{dec}
For given values of $\underline{\alpha}$ and $\rho$, according to the restrictions defined above, the space $\Sigma_{\rho,\underline{\alpha}}$ admits a unique decomposition as a union of maximal strata, say
\begin{displaymath}
\Sigma_{\rho,\underline{\alpha}}=\Sigma^{k_{1},I_{1}}\bigcup\Sigma^{k_{2},I_{2}}\bigcup\ldots\bigcup\Sigma^{k_{n},I_{n}}.
\end{displaymath}
Such a decomposition can be determined from the \textsl{S-graph} associated to $\Sigma_{\rho,\underline{\alpha}}$ and, hence,  algorithmically from the values of the parameters.
\end{fct}

A basic consequence of the selection rules given in Section 1 while defining these model spaces is the following.

\begin{fct}[\cite{cm2}]\label{inters}
Let $\Sigma^{k_{1},I_{1}}$ and $\Sigma^{k_{2},I_{2}}$ be strata that are included in $\Sigma_{\rho,\underline{\alpha}}$ for some fixed admissible values of $\underline{\alpha}$ and $\rho$. Then $\Sigma^{k_{1},I_{1}}\cap \Sigma^{k_{2},I_{2}}$ equals the union of all and only the strata that are contained both in $\Sigma^{k_{1},I_{1}}$ and $\Sigma^{k_{2},I_{2}}$, that is those $\Sigma^{k, I}$ such that $\Sigma^{k, I}\preceq\Sigma^{k_{1},I_{1}}$ and $\Sigma^{k, I}\preceq\Sigma^{k_{2},I_{2}}$.
\end{fct}

For any choice of $\left(k,I\right)$ we simply write $d_{k,I}\left(\sigma\right)=d\left(\sigma, \Sigma^{k,I}\right),$ $\sigma\in \Sigma_{\rho,\underline{\alpha}}$.
Then, for $\epsilon>0$ we define
\begin{displaymath}
\Sigma^{k,I}\left(\epsilon\right)=\left\{\sigma\in \Sigma^{k,I} | \ d_{k',I'}\left(\sigma\right)>\epsilon \ \textrm{for any couple} \ \left(k',I'\right) \ \textrm{such that} \ \Sigma^{k',I'}\prec \Sigma^{k,I} \right\}.
\end{displaymath}
\footnote{It has been pointed out by S. Galatius that the word \textsl{stratum} would probably more tipically refer to the \textsl{regular} object $\Sigma^{k,I}\left(\epsilon\right)$ rather than to $\Sigma^{k,I}$, and still we will sistematically refer here to the $\Sigma^{k,I}$s as strata for the sake of brevity.}
In case $\Sigma^{k,I}$ is such that no triplet $\left(k',I'\right)$ exists with $\Sigma^{k',I'}\prec\Sigma^{k,I}$, then we just set
\begin{displaymath}
\Sigma^{k,I}\left(\epsilon\right)=\Sigma^{k,I}.
\end{displaymath}
Such triplets $\left(k,I\right)$ will be called \textsl{minimal} with respect to $\prec$.
\newline

\begin{fct}[\cite{cm2}]\label{mnfld}
\label{mnfld}
For any couple $\left(k,I\right)$ such that the stratum $\Sigma^{k,I}$ is admissible and non-minimal and for any $\epsilon>0$ sufficiently small, the set $\Sigma^{k,I}\left(\epsilon\right)$ is a smooth open manifold of dimension $3k+card\left(I\right)-1$. As a consequence, each $\Sigma^{k,I}$ as well as the spaces $\Sigma_{\rho,\underline{\alpha}}$ is a topological finite $CW-$complex.
\end{fct}

\begin{rem}
The previous corollary involves only non-minimal strata, so one could at first wonder about minimal ones. But actually, one sees at once that they can only be of the form $\Sigma^{0,\left\{j\right\}}$ for some $j\in\left\{1,\ldots,m\right\}$. Each of these only consists of one point, so the topology of such strata is also clear.
\end{rem}

Notice that in the singular case, given $d\in\mathbb{N}$, we may have different strata having dimension $d$ whenever  $\left(k_{1},I_{1}\right)\neq\left(k_{2},I_{2}\right)$ but $3k_{1}+card\left(I_{1}\right)=3k_{2}+card\left(I_{2}\right)=d$.

\

It is also useful to recall (with a sketch of the proof) the following well-known result concerning the spaces $\Sigma^{k,\emptyset}$.

\begin{fct}[\cite{dm}]
\label{nccl}
For any $k\geq 1$, the space $\Sigma^{k,\emptyset}$ is not contractible.
\end{fct}

\begin{pfn}
We argue by induction on $k\geq 1$. When $k=1$ the space $\Sigma^{1,\emptyset}$ is trivially homeomorphic to $\Sigma$ which is assumed to be a closed 2-manifold, hence non-contractibility is clear. Let us deal with the case $k\geq 2$. In this case $\Sigma^{k,\emptyset}$ is a topological subspace of $C^{1}\left(\Sigma,g\right)$ having $\Sigma^{k-1,\emptyset}$ as topological boundary. Arguing as in \cite{bc} it is possible to prove that $\left(\Sigma^{k,\emptyset},\Sigma^{k-1,\emptyset}\right)$ is a \textsl{good pair}, and actually there is an open subset $\Omega\subset\Sigma^{k,\emptyset}$ such that:
\begin{enumerate}
\item{$\Sigma^{k-1,\emptyset}$ is a deformation retract of $\Omega$ in $\Sigma^{k,\emptyset}$;}
\item{the boundary $\partial\Omega$ is the disjoint union of $\Sigma^{k-1,\emptyset}$ and a smooth $\left(3k-4\right)-$manifold $\Lambda$.}
\end{enumerate}
Therefore $H_{\ast}\left(\Sigma^{k,\emptyset},\Sigma^{k-1,\emptyset};\Z_{2}\right)\cong \widetilde{H}_{\ast}\left(\Sigma^{k,\emptyset}/\Sigma^{k-1,\emptyset};\Z_{2}\right)$ on the one hand, $H_{3k-1}\left(\Sigma^{k,\emptyset}/\Sigma^{k-1,\emptyset};\Z_{2}\right)\cong \Z_{2}$ on the other.
We might then consider the long exact sequence in homology of the pair $\left(\Sigma^{k,\emptyset},\Sigma^{k-1,\emptyset}\right)$ and, more specifically
\begin{equation*}
\begin{CD}
\ldots @>>> H_{3k-1}(\Sigma^{k-1,\emptyset};\Z_{2})  @>>>  H_{3k-1}(\Sigma^{k,\emptyset};\Z_{2}) @>>> H_{3k-1}(\Sigma^{k,\emptyset}, \Sigma^{k-1,\emptyset};\Z_{2}) @>>> 
\end{CD}
\end{equation*}
\begin{equation*}
\begin{CD}
 @>>> H_{3k-2}(\Sigma^{k-1,\emptyset};\Z_{2}) @>>> \ldots
\end{CD}
\end{equation*}
where both $H_{3k-1}(\Sigma^{k-1,\emptyset};\Z_{2})$ and $H_{3k-2}(\Sigma^{k-1,\emptyset};\Z_{2})$ must vanish since $\Sigma^{k-1,\emptyset}$ is a finite CW-complex of dimension $3k-4$. Hence there is an isomorphism $H_{3k-1}(\Sigma^{k,\emptyset};\Z_{2})\cong H_{3k-1}(\Sigma^{k,\emptyset}, \Sigma^{k-1,\emptyset};\Z_{2})$ and so $H_{3k-1}(\Sigma^{k,\emptyset};\Z_{2})=\Z_{2}$ as claimed.
\end{pfn}

Finally, we list some properties concerning the notion of $p_{j}$-stability.

\begin{fct}\label{stab}
If $\Sigma_{\rho,\underline{\alpha}}$ is $p_{j}-$stable for some index $j$, then it is necessarily $p_{1}-$stable.
\end{fct}

\begin{pfn}
Assume $p_{j}-$stability: we need to show that given arbtrary $k, I$ the following implication holds:
\begin{displaymath}
4\pi\left[k+\sum_{i\in I}\left(1+\alpha_{i}\right)\right]<\rho \quad \Longrightarrow \quad 4\pi\left[k+\sum_{i\in \left\{1\right\}\cup I}\left(1+\alpha_{i}\right)\right]<\rho.
\end{displaymath}
There are two cases: either $j\in I$ or $j\notin I$. In the second alternative, the thesis is trivial since by assumption $\alpha_{1}\leq\alpha_{j}$ (see the paragraph on \textsl{Notations} above). In the first, define the set $\widetilde{I}$ by replacing in $I$ the index $j$ by the index $1$ (if $1\in I$, then we simply erase the index $j$). Clearly, $4\pi\left[k+\sum_{i\in\widetilde{I}}\left(1+\alpha_{i}\right)\right]<\rho$ and, thanks to the $p_{j}-$stability assumption we get $4\pi\left[k+\sum_{i\in{\left\{j\right\}}\cup \widetilde{I}}\left(1+\alpha_{i}\right)\right]<\rho$ which is equivalent to $4\pi\left[k+\sum_{i\in{\left\{1\right\}}\cup I}\left(1+\alpha_{i}\right)\right]<\rho$, so $\Sigma_{\rho,\underline{\alpha}}$ is $p_{1}-$stable.
\end{pfn}

\begin{fct}\label{ctreq}
A stratum $\Sigma^{k, I}$ is contractible if and only if the set $I$ is not emtpy and, in this case, it deformation retracts to $\delta_{p_{i}}$ for any $i\in I$. Moreover, $\Sigma^{k,I}$ is $p_{i}-$stable if and only if $i\in I$.
\end{fct}

\begin{pfn}
One implication is obvious, the other follows from Fact \ref{nccl}. The second part is a direct consequence of the algebraic version (part 2.) of Definition \ref{piuno}.
\end{pfn}

\begin{fct}\label{crit}
The space $\Sigma_{\rho,\underline{\alpha}}$ is $p_{j}-$stable for some $j$ if and only if ALL maximal strata that appear in its decomposition (see Fact 1) are $p_{j}-$stable. 
\end{fct}

\begin{pfn}
Sufficiency is trivial, while necessity follows at once by the definition of maximal stratum given above. Indeed, let $\Sigma_{\rho,\underline{\alpha}}$ be $p_{j}$-stable and let $\Sigma^{k,I}$ be one of its maximal strata. Assume by contradiction that the index $j$ does \textsl{not} belong to the set $I$. Then the stratum $\Sigma^{k,I\cup\left\{j\right\}}$ would at the same time be included in $\Sigma_{\rho,\underline{\alpha}}$ (due to $p_{j}$-stability) and \textsl{properly contain} the stratum $\Sigma^{k,I}$. But $\Sigma^{k,I}$ is maximal, contradiction.
\end{pfn}

\begin{fct}\label{pardel}
Given $k\in\mathbb{N}$, an index $j\in\left\{1,\ldots,m\right\}$ and a set $I$ such that $I\neq\emptyset$ and $j\notin I$ define for $i\in I$ the set $I(j|i)$ by replacing in $I$ the index $i$ by the index $j$. Then the topological subspace 
\begin{displaymath}
\Xi^{k,I\left(j|\cdot\right)}=\Sigma^{k,I}\bigcup\Sigma^{k,I\left(j|i_{1}\right)}\bigcup\Sigma^{k,I\left(j|i_{2}\right)}\bigcup\ldots\bigcup\Sigma^{k,I\left(j|i_{n}\right)}\subseteq C^{1}\left(\Sigma, g\right)^{\ast}
\end{displaymath}
is not contractible and, more precisely, $\widetilde{H}_{3k+card(I)-1}\left(\Xi^{k,I\left(j|\cdot\right)};\mathbb{Z}_{2}\right)=\mathbb{Z}_{2}$. 
\end{fct}

We will call for brevity such spaces $\partial\Delta$ since they formally resemble the boundary of singular simplices in the space $C^{1}\left(\Sigma, g\right)^{\ast}$.

\

\begin{pfn}
For any fixed $k\in\mathbb{N}$, we argue by induction on the cardinality of $I$. If $card(I)=1$, then $I=\left\{i\right\}$ and hence we need to show that $\widetilde{H}_{3k}\left(\Sigma^{k,\left\{i\right\}}\cup\Sigma^{k,\left\{j\right\}};\mathbb{Z}_{2}\right)=\mathbb{Z}_{2}$. This is trivial when $k=0$ (the space $\Xi^{0,I\left(j|\cdot\right)}$ being disconnected) and so assume $k\geq 1$. Observe that $\Sigma^{k,\left\{i\right\}}\cap\Sigma^{k,\left\{j\right\}}=\Sigma^{k,\emptyset}\cup\Sigma^{k-1,\left\{i,j\right\}}$ and so $H_{3k-1}\left(\Sigma^{k,\left\{i\right\}}\cap\Sigma^{k,\left\{j\right\}};\Z_{2}\right)=\Z_{2}$ as easily proved using the Mayer-Vietoris long exact sequence thanks to the previous Fact \ref{nccl} and observing that $\Sigma^{k,\emptyset}\cap\Sigma^{k-1,\left\{i,j\right\}}=\Sigma^{k-1,\left\{i\right\}}\cup\Sigma^{k-1,\left\{j\right\}}$, where $\textrm{dim}\left(\Sigma^{k,\left\{i\right\}}\right)=\textrm{dim}\left(\Sigma^{k,\left\{j\right\}}\right)=3k-3$. Again by Mayer-Vietoris, this time applied to $\Sigma^{k,\left\{i\right\}}$ and $\Sigma^{k,\left\{j\right\}}$ we get the isomorphism $H_{3k}\left(\Sigma^{k,\left\{i\right\}}\cup\Sigma^{k,\left\{j\right\}};\mathbb{Z}_{2}\right)\cong H_{3k-1}\left(\Sigma^{k,\left\{i\right\}}\cap\Sigma^{k,\left\{j\right\}};\Z_{2}\right)$ and hence the claim. Therefore, let us deal with the inductive step, $card(I)>1.$ In this case, we decompose $\Xi^{k,I\left(j|\cdot\right)}=A\cup B$ for $A=\Sigma^{k,I}$ and $B=\Sigma^{k,I\left(j|i_{1}\right)}\bigcup\Sigma^{k,I\left(j|i_{2}\right)}\bigcup\ldots\bigcup\Sigma^{k,I\left(j|i_{n}\right)}$ if $I=\left\{i_{1},\ldots,i_{n}\right\}$. Now $A$ is obviously contractible (because $I\neq\emptyset$) and $B$ is too since it is certainly $p_{j}-$stable. Moreover, 
\begin{displaymath}
A\cap B=\bigcup_{i\in I}\Sigma^{k, I\setminus\left\{i\right\}}
\end{displaymath}
and so, by inductive hypothesis (since $A\cap B$ is a $\partial\Delta$), we have $\widetilde{H}_{3k+card(I)-2}\left(A\cap B;\Z_{2}\right)=\Z_{2}$. Finally, using once again the Mayer-Vietoris long exact sequence we get the isomorpism
\begin{displaymath}
\widetilde{H}_{3k+card(I)-1}\left(\Xi^{k,I\left(j|\cdot\right)};\mathbb{Z}_{2}\right)\cong \widetilde{H}_{3k+card(I)-2}\left(\bigcup_{i\in I}\Sigma^{k, I\setminus\left\{i\right\}} ;\Z_{2}\right)
\end{displaymath} 
and this completes the proof.
\end{pfn}

\

\section{Proof of Theorem \ref{cong}}

This section is devoted to the detailed proof of Theorem \ref{cong}.

\

\begin{pfn}
We have already seen, in the Introduction, that if $\Sigma_{\rho,\underline{\alpha}}$ is $p_{1}-$stable, then it is contractible. Let us prove the converse implication. Assume then the space $\Sigma_{\rho,\underline{\alpha}}$ is \textsl{not} $p_{1}-$stable: we want to show that in this case such space is necessarily not contractible (and actually we will prove slightly more, namely that a specific homology group does not vanish). According to our structure decomposition theorem (cmp. Fact \ref{dec}) there are only two possible cases: \underline{\textsl{either}} one of the maximal strata of the decomposition is $\Sigma^{k,\emptyset}$ for some positive integer $k$ \underline{\textsl{or}} one of the maximal strata of the decomposition is $\Sigma^{k,I}$ for $I\neq\emptyset$ and $1\notin{I}$.
In order to simplify the notation we will work here with the absolute homology groups $H_{\ast}$, the argument being perfectly the same when use of reduced homology is needed.

\

\textsl{\underline{Case 1:}} Let us assume here the existence of one (and, necessarily, \textsl{only} one) maximal stratum having the form $\Sigma^{k,\emptyset}$ for some positive integer $k$. If our model space splits as
\begin{displaymath}
\Sigma_{\rho,\underline{\alpha}}=\Sigma^{k,\emptyset}\bigcup\Sigma^{k_{1},L_{1}}\bigcup\ldots\bigcup\Sigma^{k_{s},L_{s}}
\end{displaymath}
set $A=\Sigma^{k,\emptyset}$ and $B=\Sigma^{k_{1},L_{1}}\cup\ldots\cup\Sigma^{k_{s},L_{s}}$. Observe that for each index $r=1,2,\ldots,s$ we have $\textrm{dim}\left(\Sigma^{k,\emptyset}\cap \Sigma^{k_{r},L_{r}}\right)\leq 3k-3$ (as easily obtained by studying the substrata of $\Sigma^{k,\emptyset}$) and hence we might look at the following portion of the Mayer-Vietoris long exact sequence
\begin{equation*}
\begin{CD}
\ldots @>>> H_{3k-1}(A\cap B;\Z_{2})  @>>>  H_{3k-1}(A;\Z_{2})\oplus H_{3k-1}(B;\Z_{2}) @>>> H_{3k-1}(\Sigma_{\rho,\underline{\alpha}};\Z_{2}) @>>> 
\end{CD}
\end{equation*}
\begin{equation*}
\begin{CD}
 @>>> H_{3k-2}(A\cap B;\Z_{2})  @>>>  H_{3k-2}(A;\Z_{2})\oplus H_{3k-2}(B;\Z_{2}) @>>> H_{3k-2}(\Sigma_{\rho,\underline{\alpha}};\Z_{2}) @>>> \ldots
\end{CD}
\end{equation*}
which implies, due to the previous dimensional argument, that $H_{3k-1}(A)=H_{3k-1}\left(\Sigma^{k,\emptyset}\right)\cong \Z_{2}$ injects into $H_{3k-1}(\Sigma_{\rho,\underline{\alpha}};\Z_{2})$. Hence this homology group is certainly not trivial and so $\Sigma_{\rho,\underline{\alpha}}$ is not contractible, at least in this case.

\

\textsl{\underline{Case 2:}} We now turn to the case when strata $\Sigma^{k,\emptyset}$ do not appear in the decomposition of $\Sigma_{\rho,\underline{\alpha}}$, whose $p_{1}-$instability is then due to the presence of a maximal stratum $\Sigma^{k,I}$ (not necessarily only one) for $I\neq\emptyset$ and $1\notin I$. Following the construction described in proving Fact \ref{pardel}, we consider the corresponding space $\Xi^{k,I\left(1|\cdot\right)}$. Due to the very definition of our selection rules (see Definition \ref{selru}) it is seen at once that $\Xi^{k,I\left(1|\cdot\right)}\subseteq\Sigma_{\rho,\underline{\alpha}}.$ In other terms, all the strata $\Sigma^{k,I\left(1|i\right)}, i\in I$ must be strata of $\Sigma_{\rho,\underline{\alpha}}$, even though some of them might obviously be not maximal and hence these might not appear in the decomposition of such space.  We set $A=\Sigma^{k,I}$ and define $B$ as the union of all other maximal strata of $\Sigma_{\rho,\underline{\alpha}}$. Again $\textrm{dim}(A)=3k+card(I)-1$ while necessarily $H_{3k+card(I)-1}\left(A\cap B\right)=0$. Indeed, we are assuming $A$ to be a maximal stratum and therefore \textsl{none} of the strata $\Sigma^{k,I\left(1|i\right)}, i\in I$ can be contained in a stratum $\Sigma^{k+1, I\left(1|i\right)}$ of $\Sigma_{\rho,\underline{\alpha}}$, which forces in fact $A\cap B=\cup_{i\in I}\Sigma^{k, I\setminus\left\{i\right\}}$ as we studied in Fact \ref{pardel}. Looking at the Mayer-Vietoris long exact sequence, specifically
\begin{equation*}
\begin{CD}
 \ldots @>>> H_{3k+card(I)-1}(A\cap B;\Z_{2})  @>>>  H_{3k+card(I)-1}(A;\Z_{2})\oplus H_{3k+card(I)-1}(B;\Z_{2}) @>>>
\end{CD}
\end{equation*}
\begin{equation*}
\begin{CD}
  @>>> H_{3k+card(I)-1}(\Sigma_{\rho,\underline{\alpha}};\Z_{2}) @>>> H_{3k+card(I)-2}(A\cap B;\Z_{2})@>>> 
\end{CD}
\end{equation*}
\begin{equation*}
\begin{CD}
  @>>> H_{3k+card(I)-2}(A;\Z_{2})\oplus H_{3k+card(I)-2}(B;\Z_{2})@>>>\ldots
\end{CD}
\end{equation*}
we see that whenever $H_{3k+card(I)-1}(B;\Z_{2})\neq 0$ the proof is complete. Therefore, we are led to assume by contradiction $H_{3k+card(I)-1}(B;\Z_{2})=0$ and $H_{3k+card(I)-1}(\Sigma_{\rho,\underline{\alpha}};\Z_{2})=0$. In this case, the continuous inclusion $\varphi:A\cap B\hookrightarrow B$ would induce, by functoriality (and since $A=\Sigma^{k,I}$ is contractible), an \textsl{injection} $\varphi_{\ast}: H_{3k+card(I)-2}(A\cap B;\Z_{2})\rightarrow H_{3k+card(I)-2}(B;\Z_{2})$. Now, we know by Fact \ref{pardel} that $A\cap B$, being of the form $\partial\Delta$, necessarily satisfies $H_{3k+card(I)-2}(A\cap B;\Z_{2})\cong \Z_{2}\cong \left\langle \left[z\right]\right\rangle_{\Z_{2}}$ for some cycle $\left[z\right]$. On the other hand, it is checked at once that $\varphi_{\ast}\left(\left[z\right]\right)=0$ in $H_{3k+card(I)-2}(B;\Z_{2})$ (because $\left[z\right]$ is certainly trivial even in the subspace of $B$ given by $\cup_{i\in I}\Sigma^{I\left(j|i\right)}$), which contradicts injectivity of $\varphi_{\ast}$. Therefore $H_{3k+card(I)-1}(\Sigma_{\rho, \underline{\alpha}};\Z_{2})$ is not trivial and this completes the proof.
In fact, it should be remarked that such inclusion map $\varphi:A\cap B\hookrightarrow B$ is nullhomotopic and hence we can deal with a \textsl{short} exact sequence 
\begin{equation*}
\begin{CD}
 0  @>>>  H_{3k+card(I)-1}(A;\Z_{2})\oplus H_{3k+card(I)-1}(B;\Z_{2}) @>>> H_{3k+card(I)-1}(\Sigma_{\rho,\underline{\alpha}};\Z_{2})
  @>>>
\end{CD}
\end{equation*}
\begin{equation*}
\begin{CD}  
   @>>> H_{3k+card(I)-2}(A\cap B;\Z_{2}) @>>> 0
\end{CD}
\end{equation*}
so that $H_{3k+card(I)-1}(\Sigma_{\rho,\underline{\alpha}};\Z_{2})=0$ would force  $H_{3k+card(I)-2}(A\cap B;\Z_{2})=0$, contradicting Fact 8.
\
\end{pfn}
\medskip

\

\section{Proof of Theorem \ref{bh} when $m=1$ and $\rho\in\left(0,4\pi\right)$}

\medskip

\noindent

In this section we present a detailed proof of a non-existence result for equation \eqref{mod} in the special case when $m=1$ and $\rho\in\left(0,4\pi\right)$. More specifically, we prove the following:

\begin{thm}[]
\label{nonuna}
Assume that $m=1$ and $\rho\in \left(0,4\pi\right) $ is fixed. Then there exists $\alpha_{\ast}\in\left(-1,0\right)$ such that equation $\eqref{sing}_{\rho,\alpha}$ is NOT solvable for any $\alpha\leq\alpha_{\ast}$. Moreover, such $\alpha_{\ast}$ can be chosen uniformly for $\rho$ in a relatively compact subset of $\mathbb{R}_{\geq 0}\setminus 4\pi\mathbb{N}$.
\end{thm}

\
We argue by contradiction, the proof being divided into several independent steps in order to make it more readable.

\

\

\textbf{\underline{Assumption $\left(\ast\right)$:}}
suppose that the thesis is \textsl{false}, so that we can find
\begin{itemize}
\item{a sequence $\left(\alpha_{k}\right)_{k\in\mathbb{N}}$ such that $\alpha_{k}\searrow -1$:}
\item{a sequence $\left(\widetilde{u}_{k}\right)_{k\in\mathbb{N}}\subseteq C^{0}\left(\Sigma,g\right)$ such that for each $k\in\mathbb{N}$ we have that $\widetilde{u}_{k}$ is a solution of $\eqref{mod}_{\rho,\alpha_{k}}$ \textsl{and} $\int_{\Sigma}\widetilde{h}e^{2\widetilde{u}_{k}}\,dV_{g}=1$.}
\end{itemize}
Notice that the normalization condition given above can be assumed without loss of generality since our equation is invariant under traslations (i.e. $\widetilde{u}\mapsto \widetilde{u}+c$). As a consequence, each of these functions solves an equation of the form
\begin{equation}
\label{noden}
-\Delta_{g}\widetilde{u}=\rho\left(\widetilde{h}(x)e^{2\widetilde{u}}-1\right) \quad \textrm{on}\ \Sigma,
\end{equation}
which is formally simpler to handle than the previous one, \eqref{mod}.

\begin{rem}\label{indic}
\begin{enumerate}
\item{From now onwards, in order to avoid too tedious notation, we will avoid explicit use of the symbol $\widetilde{\cdot}$ in referring to solutions of the singular modified problem \eqref{mod}.}
\item{It is important to notice that the \textsl{modified datum} $\widetilde{h}$ \textsl{does depend} on the values of $\alpha_{1},\ldots,\alpha_{m}$ so that it would be more precise to adopt the notation $\widetilde{h}\left(\underline{\alpha}\right)$ and yet we will let this be implicit in order to keep the notation at a reasonable level of sophistication.}
\end{enumerate}
\end{rem}
\

\begin{lem}
\label{limita}
Let a radius $0<\d<\frac{1}{2}\textrm{inj}\left(\Sigma, p\right)$ be given. Then any sequence $\left(u_{k}\right)_{k\in \mathbb{N}}$ of normalized solutions of the problem $\eqref{mod}_{\rho,\alpha_{k}}$ for $\alpha_{k}\searrow -1$ is never uniformly bounded from below on $\partial B_{2\d}\left(p\right).$ As a consequence, any such sequence cannot be bounded in $L^{\infty}_{\textrm{loc}}\left(\Sigma\setminus \overline{B}_{\d}\left(p\right)\right).$ 
\end{lem}

\begin{pfn}
Assume by contradiction that our sequence $\left(u_{k}\right)_{k\in\mathbb{N}}$ is bounded on the boundary $\partial B_{2\d}\left(p\right),$ so that it is possible to find a constant $C\in \mathbb{R}_{>0}$ independent of $k\in\mathbb{N}$ such that 

$$
\left\{
\begin{array}{ll} \D_{g} u_{k}\leq \rho
&\text{in $ B_{2\d}\left(p\right)$;}
\\&\\
u_{k}\geq -C&\text{on $\partial B_{2\d}\left(p\right)$.}
\end{array}
\right.
$$
We are going to contradict the normalization condition $\int_{\Sigma}\widetilde{h}e^{2u_{k}}\,dV_{g}=1$ by applying a maximum/comparison principle. 

Indeed, as a first step, a standard application of the maximum (minimum) principle for the previous system implies the existence of a new constant $C\in\mathbb{R}_{>0}$ (again not depending on $k$) such that $u_{k}\geq -C$ on $B_{2\d}\left(p\right)$ for any $k\in\mathbb{N}$. As a second step, we make use of this uniform lower bound in order to apply a comparison principle. Since $\widetilde{h}(x)\simeq d_{g}(x,p)^{2\alpha}$ near  $p$ (due to standard Green's functions asymptotics, see Proposition \ref{verde}), we get the system (in local normal coordinates around the singular point $p$)

$$
\label{sistema}
\left\{
\begin{array}{ll} -\D u_{k}\geq C'\left(\left|x\right|^{2\alpha}e^{-2C}-1\right)
&\text{in $ B_{2\d}\left(p\right)$;}
\\&\\
u_{k}\geq -C&\text{on $\partial B_{2\d}\left(p\right)$.}
\end{array}
\right. 
$$

Here $C'\in\mathbb{R}_{>0}$ is a constant only depending on $\rho$, the datum $h$ and the metric $g$ (it has to be recalled here that $h$ is only assumed to be smooth and positive on $\Sigma$). Notice that in deriving the system above, we are implicitly exploiting the fact that the Laplace-Beltrami operator $\Delta_{g}$ is uniformly elliptic on $B_{\d}(p)$ and in divergence form.
As a consequence of \eqref{sistema}, a comparison argument implies the much finer estimate

\begin{displaymath}
u_{k}\left(x\right)\geq A+ \frac{C'}{4}\left(\left|x\right|^{2}-\left(2\delta\right)^{2}\right)-\frac{\l}{4\left(1+\alpha_{k}\right)^{2}}\left(\left|x\right|^{2\alpha_{k}+2}-\left(2\delta\right)^{2\alpha_{k}+2}\right)
\end{displaymath}
for some $A\in \mathbb{R}$ independent of $k$, where we have set $\l=C'e^{-2C}$. We now claim that $\int_{B_{\delta}}\widetilde{h}e^{2u_{k}}\,dV_{g}\rightarrow+\infty$ for $k\rightarrow+\infty$, which would conclude the proof of this lemma. Indeed

\begin{displaymath}
\int_{B_{2\delta}\left(p\right)}\widetilde{h}e^{2u_{k}}\,dV_{g}=2\pi \int_{0}^{2\delta}r^{1+2\alpha_{k}}e^{2A}e^{\frac{C'}{2}\left(r^{2}-\left(2\d\right)^{2}\right)-\frac{\l}{2\left(1+\alpha_{k}\right)^{2}}\left(r^{2+2\alpha_{k}}-\left(2\d\right)^{2+2\alpha_{k}}\right)}dr
\end{displaymath}

\begin{displaymath}
\geq 2\pi e^{2A}e^{{- 2C'}\d^{2}+\frac{\l}{2\left(1+\alpha_{k}\right)^{2}}\left(1-\frac{1}{2^{2+2\alpha_{k}}}\right)\left(2\d\right)^{2+2\alpha_{k}}} \int_{0}^{\d}r^{1+2\alpha_{k}}dr
\end{displaymath}
\begin{displaymath}
\geq \frac{\pi e^{2A}}{1+\alpha_{k}}{\d}^{2\left(1+\alpha_{k}\right)}e^{-2 C'\d^{2}+\frac{\l}{2}C_{\alpha_{k}}\left(2\d\right)^{2+2\alpha_{k}}}
\end{displaymath}
where $C_{\alpha_{k}}=\frac{1}{\left(1+\alpha_{k}\right)^{2}}\left(1-\frac{1}{2^{2+2\alpha_{k}}}\right)$. Therefore, once we set $\e_{k}=1+\alpha_{k}$ we get
\begin{displaymath}
\int_{B_{2\d}\left(p\right)}\widetilde{h}e^{2u_{k}}\,dV_{g}\geq \frac{C}{\e_{k}}e^{\frac{\l}{2}C_{\e_{k}}}\rightarrow +\infty
\end{displaymath}
for $k\to+\infty$ (observe that $C_{\e_{k}}\simeq \left(\log 4\right)/\e_{k} \to+\infty$, too).
The second part of our statement simply descends by observing that if we had $L^{\infty}_{\textrm{loc}}$-boundedness on $\Sigma\setminus\overline{B_{\d}}\left(p\right)$, then by standard elliptic regularity we could find a subsequence of $\left(u_{k}\right)_{k\in\mathbb{N}}$ smoothly converging to a limit function on the same domain, hence we would have boundedness from below on $\partial B_{2\d}\left(p\right)$, which is impossible by the previous argument.

\end{pfn}

To proceed further, we need to recall the following important result, that will be used to study (locally, \textsl{far} form the singularity) the behaviour of the sequence $\left(u_{k}\right)_{k\in\mathbb{N}}$ given by assumption $(\ast)$.

\begin{thm}[\cite{bm}, \cite{ls1}]\label{brezis}
Let $\Omega\subseteq\mathbb{R}^{2}$ be an open and bounded subset with smooth boundary. Consider a sequence of (classical) solutions to the equation
\begin{equation}
\label{pot}
-\Delta u_{i}=W_{i}(x)e^{2u_{i}} \quad \textrm{in}\ \Omega
\end{equation}
where $W_{i}\in C^{0}(\overline{\Omega}), \ \forall i\in\mathbb{N}$ and $0\leq W_{i}(x)\leq C_{0}$ for some fixed $C_{0}.$ Assume further that $\int_{\Omega}e^{2u_{i}}\,dx\leq C_{1}$ for some $C_{1}>0,$ not depending on the index $i.$ Then there exists a subsequence (say $(u_{i_{k}})$) which satisfies one of the following three alternatives:
\begin{itemize}
\item{$(u_{i_{k}})$ is uniformly bounded;}
\item{$(u_{i_{k}})\to-\infty$ uniformly on any compact subset of $\Omega;$}
\item{there exists a finite blow-up set $S=\left\{a_{1},...,a_{m}\right\}$ such that for any $j\in\left\{1,...,m\right\}$ there exists a sequence $\left(x_{j_{k}}\right)_{k\in\mathbb{N}}\subseteq\Omega$ converging to $a_{j}$ with $u_{i_{k}}(x_{j_{k}})\to+\infty$ but for any compact subset of $\Omega\setminus S$ the functions $u_{i_{k}}$ converge uniformly to $-\infty.$ Moreover, in this case, there exist constants $\beta_{1},...,\beta_{m}$ such that $\beta_{j}\geq 4\pi$ for any $j\in\left\{1,\ldots,m\right\}$ and $W_{i_{k}}e^{2u_{i_{k}}}\rightharpoonup \sum_{j=1}^{m}\beta_{j}\delta_{a_{j}}$ in the sense of measures.}
\end{itemize}
\end{thm}

This is indeed the crucial tool in the next step of our proof. 

\
 
In the following discussion, it will be often convenient to write our equation in \textsl{isothermal coordinates} and perform a suitable change of variables that we introduce here. Given  any $q\in\Sigma$ and a small ball $B=B_{r}(q)\subset\Sigma$ (for our purposes we can certainly assume $r<\frac{1}{2} inj\left(\Sigma, q\right)$) it is well-known (see for instance \cite{sp}) that there exists a system of coordinates $z=z(x): B_{2r}(q) \to \Omega\subset{\mathbb{R}^{2}}$ around the point $q$ and a function (called \textsl{\textsl{conformal factor}}) $w\in C^{\infty}\left(\Omega\right)$ such that $w(0)=0$ and $g(z)=e^{-2w}\left(dz^{1}\otimes dz^{1}+ dz^{2}\otimes dz^{2}\right)$. Therefore, the (locally defined) Riemannian metric $\widehat{g}=e^{2w}g$ is \textsl{flat} on $B_{r}\left(q\right)$. Clearly, we can extend such $w$ to a smooth function on the whole $\mathbb{R}^{2}$ which vanishes identically on $\mathbb{R}^{2}\setminus\Omega$. Let then $T\in C^{\infty}\left(\mathbb{R}^{2}\right)$ a solution of the problem $\D T=\rho e^{-2w}$ such that $T(0)=0$. Making use of the fact that $\D_{\widehat{g}}=e^{-2w}g$ we immediately get that if $u$ solves \eqref{noden}, then the function $v=u-T$ satisfies on $B$ the equation
\begin{equation}
\label{nuova}
-\D v=\rho V e^{2v}
\end{equation}
where $V=\widetilde{h}e^{-2w+2T}$.

\begin{lem}
\label{altern}
Assume the hypotheses $(\ast)$ and let a radius $0<\d<\frac{1}{2}\textrm{inj}\left(\Sigma, p\right)$ be given. Then the following alternative holds: \newline
\

\underline{\textsl{either}} there exists a subsequence $\left(u_{k_{l}}\right)_{l\in\mathbb{N}}$ of $\left(u_{k}\right)_{k\in\mathbb{N}}$ and a smooth function $u_{\infty}\in C^{\infty}\left(U\right)$ for some open set $U\supseteq\Sigma\setminus B_{\d}\left(p\right)$ such that $u_{k_{l}}\rightarrow u_{\infty}$ uniformly on $\Sigma\setminus B_{\d}\left(p\right)$; \newline
\

\underline{\textsl{or}} there exists a subsequence $\left(u_{k_{l}}\right)_{l\in\mathbb{N}}$ such that $u_{k_{l}}\rightarrow -\infty$ uniformly on $\Sigma\setminus B_{\d}\left(p\right)$.
\end{lem}

\

\begin{pfn}
By a standard compactness argument, we can always find a finite number of (geodesic) balls of $\Sigma$, say $B_{1},\ldots, B_{N}$ such that $\Sigma\setminus B_{\d}\left(p\right)\subseteq\cup_{i=1}^{N} B_{i}$ and $\cup_{i=1}^{N} 2B_{i}\subseteq \Sigma\setminus B_{\d/2}\left(p\right)$. Moreover, arguing as above, we can get on each $2 B_{i}$ a smooth function $w_{i}$ so that the conformal metric $\widehat{g}=e^{2w_{i}}g$ is \textsl{flat}. We remark here that this can not be done globally, namely on the whole $\Sigma\setminus B_{\d}\left(p\right)$ because the \textsl{Uniformization Theorem} only allows solvability of the Yamabe problem on $\Sigma$ (that is $K_{\widehat{g}}=const.$ on $\Sigma$). However, in such local isothermal coordinates $z_{i}=z_{i}(x): 2B_{i}\to\Omega_{i}\subset\mathbb{R}^{2}$ we get the equation
\begin{displaymath}
-\D v_{k} =\rho V_{k}^{i}e^{2v_{k}} \quad \textrm{on}\ \Omega_{i},
\end{displaymath}
by solving the Dirichlet problems $\D T_{i}=\rho e^{-2w_{i}}$ and setting $v_{k}=u_{k}-T_{i}$ on $\Omega_{i}$, $V_{k}^{i}= e^{2z_{i}}\widetilde{h}$ where clearly $z_{i}=-w_{i}+T_{i}$.
Now, all of the data $V_{k}^{i}$ are both smooth, positive and uniformly bounded respectively on each $2B_{i}$. Therefore, we are then in position to apply the blow-up analysis by Brezis and Merle.
Notice that, due to the fact that $\rho\in\left(0,4\pi\right)$ and that $\int_{\Omega_{i}}V_{k}^{i}e^{2v_{k}}\,d\mathscr{L}^{2}\left(z\right)\leq 1$ (as immediately checked from the normalization condition $\left(\ast\right)$), \textsl{the third alternative can be ruled out} and therefore, by some standard elliptic regularity argument, we are led to the following alternative:\newline
\underline{\textsl{either}} there exists a subsequence $\left(u_{k_{l}}\right)_{l\in\mathbb{N}}$ and a smooth function $u_{\infty}^{i}\in C^{\infty}\left(2B_{i}\right)$ such that $u_{k_{l}}\rightarrow u_{\infty}^{i}$ uniformly on $\overline{B_{i}}$; \newline
\underline{\textsl{or}} there exists a subsequence $\left(u_{k_{l}}\right)_{l\in\mathbb{N}}$ such that $u_{k_{l}}\rightarrow -\infty$ uniformly on $\overline{B_{i}}$.

Obviously, for each $B_{i}$ there exists a $B_{i'}$ such that $B_{i}\cap B_{i'}\neq\emptyset$ snd hence we can certainly label the balls so that for each $i=2,\dots, N$ there is a \textsl{smaller} index $i'$ such that $B_{i}\cap B_{i'}\neq\emptyset$: as a consequence, due to incompatibility of the previous two cases, a (finite) diagonal argument allows to obtain the previous alternative \textsl{globally} on $\Sigma\setminus B_{\delta}(p)$. This precisely means that we can find a subsequence of $\left(u_{k}\right)_{k\in\mathbb{N}}$ (that for the sake of clarity is not renamed) such that one and only one of the alternatives in our statement holds and so the proof is complete.

\end{pfn}

\

Hence, if we put together Lemma \ref{altern} and Lemma \ref{limita}, we get that the asymptotic behaviour for $k\to+\infty$ of our sequence is \textsl{uniquely determined}.

\begin{lem}
\label{unica}
Assuming $\left(\ast\right)$, we have that for any $0<\d<\frac{1}{2}\textrm{inj}\left(\Sigma,p\right)$ there exists a subsequence of $\left(u_{k}\right)_{k\in\mathbb{N}}$, say $\left(u_{k_{l}}\right)_{l\in\mathbb{N}}$,  such that 
\begin{displaymath}
u_{k_{l}} \longrightarrow -\infty \quad \textrm{uniformly} \ \textrm{on} \ \Sigma\setminus B_{\d}\left(p\right) \textrm{for} \ l\to +\infty.
\end{displaymath}
\end{lem}

As a further consequence, by applying a standard diagonal argument to $\left(u_{k}\right)_{k\in\mathbb{N}}$ for a \textsl{discrete} sequence of values of $\d$ monotonically converging to $0$ we can prove the following easy, but enligthning, result.

\begin{lem}
\label{masse}
There exists a subsequence of $\left(u_{k}\right)_{k\in\mathbb{N}}$ (which we do not rename) such that $\mu_{k}\rightharpoonup \delta_{p}$ (weakly, in the sense of measures), where we have set $\mu_{k}=\widetilde{h}e^{2u_{k}}\,dV_{g}$.
\end{lem} 

\begin{pfn}
As explained above, it is possible by means of the previous lemmas to get a subsequence such that the associated (probablity) measures $\left(\mu_{k}\right)_{k\in\mathbb{N}}$ satisfy this property:
\begin{displaymath}
\textrm{for each compact set} \ \Xi\subseteq\Sigma\setminus{\left\{p\right\}}, \quad \mu_{k}\left(\Xi\right)\to 0.
\end{displaymath}
Hence, the family $\left(\mu_{k}\right)_{k\in\mathbb{N}}\subseteq\mathcal{P}(\Sigma)$ is \textsl{tight} and therefore, by Prohorov's compactness theorem, there exists a Borel probability measure $\mu$ such that $\mu_{k}\rightharpoonup\mu$ (up to possibly extracting an extra subsequence). Clearly, such a measure $\mu$ can only be supported at $p$ and so the conclusion follows.
\end{pfn}

The core of our argument relies on Green's representation formula, which is recalled below.

\begin{pro}{(\cite{aub})}
\label{verde}
Let $\left(M,g\right)$ a closed, smooth Riemannian manifold of dimenion 2 and having volume $V$. The Green function $G(x,y)\in C^{\infty}\left(M\times M\setminus \textrm{diag}_{M\times M}\right)$ for the Laplace operator $\Delta_{g}$ is defined (up to an additive costant) as the solution of the problem
\begin{displaymath}
\Delta G(x,y)=\d_{x}(y)-V^{-1}.
\end{displaymath}
and possesses the following properties:
\begin{enumerate}
\item{For all functions $\varphi\in C^{2}(M,g)$
\begin{displaymath}
\varphi(x)=V^{-1}\int_{M}\varphi(y)\,dV_{g}(y)+\int_{M}G(x,y)\D_{g}\varphi(y)\,dV_{g}(y)
\end{displaymath}
and this equality holds whenever such integrals make sense;}
\item{There exists a constant $C$ (depending on $M$ and the metric $g$) such that
\begin{displaymath}
\left|G(x,y)\right|<C\left(1+\left|\log r\right|\right), \quad \left|\nabla_{g}G(x,y)\right|<Cr^{-1}
\end{displaymath}
where $r$ denotes the Riemannian distance from the point $x\in M$, namely $r=d_{g}\left(x,y\right)$;}
\item{The value of the integral $\int_{M}G(x,y)\,dV_{g}\left(y\right)$ does not depend on $x$ and hence we can choose Green's function so that its integral is zero.}
\item{If $d_{g}$ is the Riemannian distance on the 2-manifold $M$, then $G(x,y)\simeq \log d_{g}\left(x,y\right)$ for $x\to y$.}
\end{enumerate}
\end{pro}

In the sequel, we always implicitly assume that the definition of the Green function $G(x,y)$ is given with the additional requirement that $\int_{M}G(x,y)\,dV_{g}\left(y\right)=0$. As a consequence of this classic result, we can represent each function obtained by means of the diagonal argument explained before Lemma \ref{masse} as follows:
\begin{equation}
\label{intermedia}
u_{k}\left(x\right)= \overline{u}_{k}-\rho\int_{\Sigma}G\left(x,y\right)\left(\widetilde{h}\left(y\right)e^{2u_{k}}\left(y\right)\right)\,dV_{g}\left(y\right)
\end{equation}
and the idea now is to observe that the convergence result Lemma \ref{masse} implies that $u_{k}-\overline{u}_{k}$ \textsl{looks like} the Green function $G(p,\cdot)$, at least \textsl{far} from the singularity. Let us see how this can be formalized.

\begin{lem}
\label{conti}
Let $\left(u_{k}\right)_{k\in\mathbb{N}}$ the sequence obtained in Lemma \ref{masse} assuming $\left(\ast\right)$ and let $0<\eta<\frac{1}{2}\textrm{inj}\left(\Sigma, p\right)$. Then the following representation formula holds for each function $u_{k}$:
\begin{displaymath}
u_{k}\left(x\right)=\overline{u}_{k}-\rho G\left(x,p\right)+R_{1,k}\left(x\right)
\end{displaymath}
where 
\begin{displaymath}
R_{1,k}, \ \nabla_{g}R_{1,k} \ \longrightarrow 0 \ \textrm{uniformly on the annulus} \ S=\left\{\frac{\eta}{2}\leq d_{g}\left(p,\cdot\right)\leq\frac{3\eta}{2}\right\} \ \textrm{for} \  k\to +\infty.
\end{displaymath}
\end{lem}

As indicated by the referee, the proof of this Lemma is quite standard given the representation formula \eqref{intermedia} and the result contained in Lemma \ref{masse}, so we finally decided to omit it. In fact, it is gotten by means of elliptic estimates to be used together with Proposition \ref{verde}.

\begin{rem}
It has to be noticed that Green's functions are extensively used in Geometric Analysis, the most relevant example being the solution of the Yamabe problem by Schoen and Yau by means of the Positive Mass Theorem. In the study of mean-field equations several applications might be mentioned, such as for instance \cite{bm} and \cite{clin}.
\end{rem}

The previous lemma allows us to estimate the gradient and the normal derivative $\nabla v, \left(\nu,\nabla v\right)$ for $z\in\partial B_{\d}\left(p\equiv 0\right)$ in one of the \textsl{leading} terms of the Pohozaev identity which is stated below. This is the crucial issue in finding a contradiction from assumption $\left(\ast\right)$, which concludes our argument.

\

\begin{lem}\label{poho}
Suppose that $v$ is a solution of equation \eqref{nuova}, let $q\in\Sigma$ and let $\eta< \frac{1}{2}\textrm{inj}\left(\Sigma, q\right)$, then the following identity holds:
\begin{equation}
\label{pohol}
\int_{\left|z\right|=\d}\left(z,\nu\right)\rho Ve^{2v}\,d\sigma\left(z\right)+\int_{\left|z\right|=\d}\left(z,\nu\right)\left[2\left(\nu,\nabla v\right)^{2}-\left|\nabla v\right|^{2}\right]\,d\sigma\left(z\right)
\end{equation}
\begin{equation}
\label{pohor}
=\rho\left[\int_{\left|z\right|<\d}2Ve^{2v}d\mathscr{L}^{2}\left(z\right)+\int_{\left|z\right|<\d}\left(z,\nabla \log V\right)Ve^{2v}d\mathscr{L}^{2}\left(z\right)\right]
\end{equation}
where $z: B_{\eta}\left(q\right)\to\Omega\subset\mathbb{R}^{2}$ are isothermal coordinates (see the discussion before equation \eqref{nuova}), $\d=\d(\eta)$ is a suitably small radius such that $B_{\d}\left(0\right)\subset\Omega$, $\left(\cdot,\cdot\right)$ is the dot product in $\mathbb{R}^{2}$ and $\nu$ is the outer unit normal.
\end{lem}

\

Thanks to the preliminary reduction to the Euclidean setting, the proof of this result is quite standard and essentially follows the argument presented for instance in Chapter 8 of \cite{am}. 

\

\begin{pfn}(Theorem \ref{nonuna}).
Assuming $\left(\ast\right)$, we have already seen through all the previous steps that we can always find a sequence $\left(u_{k}\right)_{k\in\mathbb{N}}$ solving \eqref{mod} for $\alpha_{k}\searrow -1$ and satisfying the asymptotic characterization given by Theorem \ref{conti}. Now, fix an appropriate, \textsl{small} radius $0<\eta<\frac{1}{2}\textrm{inj}\left(\Sigma, p\right)$ and  let us work in isothermal coordinates $z=z(x): B_{\eta}(p)\to\Omega\subset\mathbb{R}^{2}$ (around the singular point $p=p_{1}$) and adopt the very same notations introduced before equation \eqref{nuova}. So we get a sequence of functions $\left(v_{k}\right)_{k\in\mathbb{N}}$ and potentials $V_{k}$ defined on $\Omega$ and such that for any index $k\in\mathbb{N}$ the following equation holds:
\begin{displaymath}
\label{mitica}
-\D v_{k}=\rho V_{k} e^{2v_{k}}, \quad z\in\Omega.
\end{displaymath}
By using these coordinates and the definition of the Green's function $G(x,p)$ (Proposition \ref{verde}), it is straightforward to check that its local form $G_{0}(z,0)$ reads
\begin{displaymath}
G_{0}(z,0)=\frac{1}{2\pi}\log \left|z\right|+\sigma(z), \quad  z\in B_{2\d}\left(0\right)\subset\Omega
\end{displaymath}
for some smooth function $\sigma$. 
Then, using Lemma \ref{masse}, we get that
\begin{equation}
\label{cambia}
V_{k}e^{2v_{k}}\rightharpoonup \d_{z=0}
\end{equation}
weakly in the sense of measures in $\Omega$ (and uniformly on any compact subset of $\Omega\setminus\left\{0\right\}$), while using the asymptotic Lemma \ref{conti} we obtain
\begin{equation}
\label{deriv}
\n v_{k}\left(z\right)=-\frac{\rho}{2\pi}\frac{z}{\left|z\right|^{2}}+R_{2,k}\left(z\right), \quad R_{2,k}\left(z\right)\to R_{2}\left(z\right) \ (\textrm{for } \ k\to+\infty),
\end{equation}
the convergence being uniform in $\overline{B}_{2\eta}\left(0\right)$ with $R_{2}$ (say) continuous. The Pohozaev identity for $v_{k}$  on $B_{\eta}\left(0\right)$ (in the form recalled above, as Lemma \ref{poho}) reads
\begin{displaymath}
\int_{\left|z\right|=\d}\left(z,\nu\right)\rho V_{k}e^{2v_{k}}\,d\sigma\left(z\right)+\int_{\left|z\right|=\d}\left(z,\nu\right)\left[2\left(\nu,\nabla v_{k}\right)^{2}-\left|\nabla v_{k}\right|^{2}\right]\,d\sigma\left(z\right)=
\end{displaymath}
\begin{displaymath}
\rho\left[\int_{\left|z\right|<\d}2V_{k}e^{2v_{k}}d\mathscr{L}^{2}\left(z\right)+\int_{\left|z\right|<\d}\left(z,\nabla \log V_{k}\right)V_{k}e^{2v_{k}}d\mathscr{L}^{2}\left(z\right)\right].
\end{displaymath}
At this point, by using \eqref{cambia}, \eqref{deriv} and the fact that $\log V_{k}=2\alpha \log\left|z\right|+\psi_{k}(z)$ for some smooth functions $\psi_{k}$ (these being $C^{1}-$uniformly bounded on $B_{2\d}(0)$) we can easily evaluate the leading terms in this identity and get
\begin{displaymath}
\frac{1}{2\pi}\rho^{2}=2\left(1+\alpha_{k}\right)\rho+ R_{3,k}\left(\d\right),
\end{displaymath}
where $\lim_{\eta\to 0^{+}}\lim_{k\to+\infty} R_{3,k}\left(\d\right)=0$. 
Clearly, as we let $\alpha_{k}\searrow -1$ we obtain $\rho=0$ which is desired contradiction. 

\end{pfn}

\begin{rem}
\begin{itemize}
\item{It should be highlighted that it is possible to give a proof of this theorem (\ref{nonuna}) without making use isothermal coordinates. In fact, in our first proof of this result all these computations were performed in normal coordinates around the singular point $p$. The drawback of such a strategy is that an \textsl{ad hoc} local Pohozaev identity is required in that case.} 
\item{It was pointed out to us that these sorts of arguments, based on the analysis of the Pohozaev identity near the concentration point were introduced in \cite{bt1}, see also \cite{bclt} and \cite{bmt}.}
\end{itemize}
\end{rem}

\medskip

\

\section{Proof of Theorem \ref{bh} in the general case}

In this section, we describe how the previous strategy needs to be modified in order to obtain the general non-existence result stated in Section 1. We remark, as discussed in the paragraph concerning \textsl{Notations} that it is possible to assume, without loss of generality, that $j=1$.

\medskip

\begin{pfn}
First of all, observe that Lemma \ref{limita} applies without changes to this more general framework. Hence, we may then argue by a local covering argument as in Lemma \ref{altern} (with some balls covering the singularities $p_{2},\ldots,p_{m}$) and thanks to an appropriate extension of Theorem \ref{brezis} to the singular setting (given by Theorem 2.1 in \cite{bmt}), we finally end up with the following (global) alternative on $\Sigma\setminus B_{\d}\left(p\right)$ (for each $\d>0$):\newline
\underline{\textsl{either}} there exists a subsequence $\left(u_{k_{l}}\right)_{l\in\mathbb{N}}$ of $\left(u_{k}\right)_{k\in\mathbb{N}}$ that blows-up at finitely many points $a_{1},\ldots, a_{N}$ in the sense described above; \newline
\underline{\textsl{or}} there exists a subsequence $\left(u_{k_{l}}\right)_{l\in\mathbb{N}}$ such that $u_{k_{l}}\rightarrow -\infty$ uniformly on $\Sigma\setminus B_{\d}\left(p\right)$.

Now, these two cases are not essentially distinct since the second may be just recovered by setting formally $N=0$ in the first.
Anyhow, we are in position to apply a diagonal argument as in Lemma \ref{masse} and get
\begin{displaymath}
\rho \mu_{k}\rightharpoonup 4\pi\sum_{i=1}^{N}\o_{i}\d_{a_{i}}+M\d_{p}
\end{displaymath}
where we recall that the blow-up analysis (see \cite{ls1}, \cite{yanyan1}, \cite{bt1}, \cite{bmt}) forces \textsl{quantization of mass}, namely
\begin{displaymath}
\o_{i}=\left\{\begin{matrix} 1+\alpha_{j} & \mbox{if }a_{i}=p_{j} \ \mbox{for some} \ j=2,\ldots,m; \\ 1 & \mbox{otherwise}.
\end{matrix}\right.
\end{displaymath}
and we have set
\begin{displaymath}
M=\rho-4\pi\sum_{i=1}^{N}\o_{i}.
\end{displaymath}
The trivial but crucial remark here is that
\begin{displaymath}
\rho\notin \mathfrak{S}\left(\alpha_{2},\ldots,\alpha_{m}\right) \ \Rightarrow \ M>0
\end{displaymath}
and this remark will be useful to get a contradiction by applying Pohozaev inequality.
The next step is to apply Green's representation formula for each function $u_{k}$. Let us fix a strip $S$ (similarly to what we did in Lemma \ref{conti}) around the point $p=p_{1}$ and far from the other singularities. As a consequence we can represent, for $x\in S$, the sequence $\left(u_{k}\right)_{k\in\mathbb{N}}$ as follows:
\begin{equation}
\label{resti}
u_{k}\left(x\right)=\overline{u}_{k}-4\pi \sum_{i=1}^{N}\o_{i}G\left(x,a_{i}\right)-M G\left(x,p\right)+R_{4,k}\left(x\right)
\end{equation}
where 
\begin{displaymath}
R_{4,k}, \ \nabla_{g}R_{4,k} \ \longrightarrow 0 \ \textrm{uniformly on the strip} \ S=\left\{\frac{\eta}{2}\leq d_{g}\left(p,\cdot\right)\leq\frac{3\eta}{2}\right\} \ \textrm{for} \  k\to +\infty.
\end{displaymath}
Again, notice that the remainder terms depend on the choice of $\eta$, as above.
After the introduction of local isothermal coordinates around the singular point $p_{1}$, we are finally in position to apply the Pohozaev identity, Lemma \ref{poho}, to our problem in the form \eqref{nuova}. 
To this aim, observe that in the representation formula \ref{resti} the term $\sum_{i=1}^{N}\o_{i}G\left(x,a_{i}\right)$ is (say) $C^{1}-$uniformly bounded on $S=S(\eta)$ (and for all sufficiently small $\eta$), so that we are led to the very same situation treated above in the proof of Theorem \ref{nonuna} and we can write
\begin{equation}
\label{semplice}
\nabla v_{k}\left(z\right)=-\frac{M}{2\pi}\frac{z}{\left|z\right|^{2}}+R_{5,k}\left(z\right), \quad R_{5,k}\left(z\right)\to R_{5}\left(z\right) \ (\textrm{for } \ k\to+\infty).
\end{equation}
As a consequence, arguing as above we end up getting
\begin{displaymath}
\frac{M^{2}}{2\pi}=2\left(1+\alpha_{k}\right)M+R_{6,k}\left(\d\right)
\end{displaymath}
for a suitable remainder term satisfying $\lim_{\d\to 0^{+}}\lim_{k\to+\infty}R_{6,k}\left(\d\right)$.
Therefore letting first $k\to+\infty$ and then $\d\to 0$ we obtain $M=0$, which contradicts the assumption $\rho\notin \mathfrak{S}\left(\alpha_{2},\ldots,\alpha_{m}\right)$.
\end{pfn}

\medskip

\medskip

\

\

\textsl{Acknowledgments}. The author would like to thank A. Malchiodi for constant support, encouragement and many fruitful discussions concerning the Liouville equation in presence of conical singularities. Some relevant ideas concerning the non-existence argument presented in Section 4 are due to him. Gratitude is also expressed to S. Galatius for his helpful suggestions and for carefully reading the manuscript. We are indebted to the referee for several improvements and more specifically for suggesting a simplified proof of Theorem \ref{bh}. During the preparation of this work, the author was supported by the FIRB-Ideas project {\em Analysis and Beyond} from
MiUR and by NSF grant DMS/0604960.

\end{document}